\documentclass[12pt,a4paper]{amsart}

\usepackage{bbm}

\usepackage
{hyperref}

\newcommand {\fbgl}{{\mathfrak{bgl}}}

\newcommand {\fwk}{{\mathfrak{wk}}}

\newcommand {\htr}{{\text{htr}}}

\newcommand{\Size}{\text{Size}}

\def\mmat #1,#2,#3,#4,{\text{\small\arraycolsep=3pt $
\begin{pmatrix}#1&#2\\#3&#4\end{pmatrix}$}}

\usepackage{DLdef1}
\let\ssec\subsection

\renewcommand {\ssbegin}[2][*]
 {\refstepcounter{subsection}%
\if#1*
\addcontentsline{toc}{subsection}{\thesubsection.\hskip 1pc #2}%
\else
\addcontentsline{toc}{subsection}{\thesubsection.\hskip 1pc #2. #1}%
\fi
 \def \secno {\gdef \secno {}{\ssecfont
\thesubsection.\hskip 2ex}%
 }%
 \begin{#2}}

\renewcommand {\sssbegin}[2][*]
  {\refstepcounter{subsubsection}
\if#1*
\addcontentsline{toc}{subsubsection}{\thesubsubsection.\hskip 1pc #2}%
\else
\addcontentsline{toc}{subsubsection}{\thesubsubsection.\hskip 1pc #2. #1}
\fi
  \def \secno {\gdef \secno {}{\ssecfont \thesubsubsection.\hskip 2ex}%
  }%
   \begin{#2}}

\renewcommand {\parbegin}[2][*]
  {\refstepcounter{paragraph}
\if#1*
\addcontentsline{toc}{paragraph}{\theparagraph.\hskip 1pc #2}%
\else
\addcontentsline{toc}{paragraph}{\theparagraph.\hskip 1pc #2. #1}
\fi
  \def \secno {\gdef \secno {}{\ssecfont \theparagraph.\hskip 2ex}%
  }%
   \begin{#2}}

\begin{document}

\title[NIS in characteristic 2]{Nondegenerate invariant symmetric bilinear forms on simple Lie superalgebras in characteristic 2}

\author[A.~Krutov]{Andrey Krutov}
\address{Institute of Mathematics, Czech Academy of Sciences, \v{Z}itn\'a 25, 115 67 Prague, Czech Republic}
\email{krutov@math.cas.cz}

\author[A.~Lebedev]{Alexei Lebedev}
\address{Equa Simulation AB, R{\aa}sun\-da\-v\"agen 100, Solna, Sweden}
\email{alexeylalexeyl@mail.ru}

\author[D.~Leites]{Dimitry Leites}
\address{New York University Abu Dhabi, Division of Science and Mathematics, P.O. Box 129188, United Arab
  Emirates}
\address{Department of mathematics, Stockholm University, SE-106 91 Stockholm, Sweden}
\email{dl146@nyu.edu, mleites@math.su.se}

\author[I.~Shchepochkina]{Irina Shchepochkina}
\address{Independent University of Moscow, Bolshoj Vlasievsky per, dom
  11, RU-119 002 Mos\-cow, Russia}
\email{irina@mccme.ru}

\thanks{A.K. was supported through the program ``Oberwolfach Leibniz Fellows'' by
the Mathematisches Forschungsinstitut Oberwolfach (MFO) in 2018, by the QuantiXLie Centre of Excellence,
a~project cofinanced by the Croatian Government and European Union through the European Regional Development
Fund --- the Competitiveness and Cohesion Operational Programme (KK.01.1.1.01.0004), and by the GA\v{C}R
project 20-17488Y and RVO: 67985840. D.L. was supported by the grant AD 065 NYUAD.
A.K. and D.L. heartily thank MFO for hospitality and excellent working conditions.}



\begin{abstract} As is well-known, the dimension of the space spanned by the non-degenerate invariant symmetric bilinear forms (NISes) on any simple finite-dimensional Lie algebra or Lie superalgebra is equal to at most 1 if the characteristic of the algebraically closed ground field is not~2.

We prove that in characteristic 2, the superdimension of the space spanned by NISes can be equal to 0, or 1, or $0|1$, or $1|1$;
it is equal to $1|1$ if and only if the Lie superalgebra is a~queerification (defined in arXiv:1407.1695) of a~simple classically restricted Lie algebra with a~NIS (for examples, mainly in characteristic $\neq 2$, 
see arXiv:1806.05505).

We give examples of NISes on deformations (with both even and odd  parameters) of several simple finite-dimensional Lie superalgebras in characteristic 2.

We also recall examples of multiple NISes on simple Lie algebras over non-closed fields.
\end{abstract}

\maketitle

\thispagestyle{empty}

\section{Introduction}

This paper is a~sequel to \cite{BKLS}, where non-degenerate invariant symmetric bilinear forms (briefly: NISes) on   known simple  Lie algebras and Lie superalgebras (finite-dimensional and $\Zee$-graded of polynomial growth) are listed if the characteristic $p$ of the ground field $\Kee$ is distinct from 2, and for occasional examples where $p=2$. Here we mainly consider an algebraically closed field $\Kee$ for $p=2$. For numerous applications of Lie (super)algebras with a~NIS, see \cite{DSB, BeBou,BLS}.
In particular, NISes allow one to distinguish true deforms, see for example~\cite[Claims~3.3 and~3.4]{BKLS};
this is crucial in the classification of simple finite-dimensional Lie (super)algebras over algebraically
closed fields of characteristics~2, see~\cite{BLLS}.
In characteristics~0, the existence of NISes on a (not necessary simple) Lie
algebra~$\fg$ allows one to recover a~lot of information about the structure of~$\fg$; for examples,
see~\cite{A,DPU} and references therein. The dimension of the space of even
NISes on a~Lie super algebra is an~important invariant; see~\cite{Be2}.
We are thankful to the referee who pointed at the following fact shown in~\cite{Be1}: if a Lie superalgebra $\fg=\fg_\ev\oplus\fg_\od$ with an even NIS is
such that the odd part~$\fg_\od$ is a~completely reducible module over a non-zero even part~$\fg_\ev$, then
$\fg$ is simple if and only if such even NIS in unique (up to multiplication by a scalar).
Applications of odd NISes were considered in~\cite{ABB}.

Our main result is a~newly discovered fact in characteristic 2 ---  the description of the possible
superdimension of the space of NISes on a~given simple finite-dimensional Lie superalgebra,  see Theorem~\ref{th}. For descriptions of simple finite-dimensional Lie superalgebras over algebraically closed field of characteristic 2,  see \cite{BGLLS} and references therein.

For completeness, we consider Theorem \ref{th} for any characteristic since we do not know if it was ever
published in full generality, although known as a~folklore for $p\neq2$. If $p\ne 2$, Theorem~\ref{th}  seems to be \lq\lq
just a~direct generalization" of a~well-known fact about Lie algebras,
\textit{unless we realize that we might encounter objects  \textbf{inhomogeneous} with respect to the parity}.  For $p=2$, the situation is even more delicate; its complete description is a~new result, cf.~\cite{SF}.

For results on NISes on simple finite-dimensional Lie algebras over algebraically closed fields of characteristic $p\neq 2$,  see~\cite{BZ, Bl, GP} and~\cite{BKLS}.

\ssec{Generalities}
From the super $K_0$-functor point of view, see \cite{Mi}, the \textit{superdimension} of a~given superspace $V$ is $\sdim V:=\dim V_\ev+\eps\dim V_\od$, where $\eps^2=1$; usually one writes $\sdim V=\dim V_\ev\,|\,\dim V_\od$, so $\eps=0|1$. We write $\dim V:=\dim V_\ev+\dim V_\od$.

Let $\fg$ be a~Lie superalgebra over~$\Kee$. A bilinear form~$B$ on~$\fg$ is called \emph{homogenous} (with
respect to parity) if

1) $\fg_\ev \perp \fg_\od$, in this case $B$ is called \emph{even};

2) $\fg_\ev$ and $\fg_\od$ are isotropic subspaces, in this case $B$ is called \emph{odd}.

Recall that a~bilinear form $B$ is called \emph{invariant} if the following condition holds
\[
  B([x,z],y) - B(x,[z,y]) = 0\qquad\text{for all $x,y,z\in\fg$}.
\]

\ssec{Conventions in this note} 1) Speaking about ``the space of NISes'', we exercise the usual abuse of the language: we are speaking about the space spanned by all NISes, but not all forms in this space have to be non-degenerate. For example, the zero form is never non-degenerate. Observe that all forms in the space spanned by all NISes are, however, invariant and symmetric; they are briefly called ISes.

In the case where ``the space of NISes'' is of superdimension $1|1$, we mean that, up to a~nonzero factor, there is one even non-degenerate form and one odd non-degenerate form. This space of NISes contains a~1-dimensional inhomogeneous (with respect to parity) space of degenerate forms (ISes).

2) The ground field is algebraically closed  of characteristic $p$ (mainly, $p=2$ or 0), unless otherwise specified.

3) We consider only finite-dimensional (super)algebras.

4) All (super)commutative (super)algebras are supposed to be associative with~1; their morphisms  should  send 1 to 1, and  the morphisms  of supercommutative superalgebras should preserve parity.

\subsection{On the contents of this paper}\label{cont}
Section \ref{sec:main} contains the Main Theorem~\ref{th} --- the most interesting result of the paper.

Section \ref{basics}
contains vital material: even the definition of Lie superalgebra for $p=2$ (and 3) is not obvious.
We also recall definition of classical restrictedness of Lie superalgebra for $p=2$; for non-classical ones, see \cite{BLLS}.

Section \ref{sec:examples}
contains several examples illustrating the Main Theorem~\ref{th}.

In Section\ref{sec:remarks}, we examine existence of multiple NISes over non-closed fields and define a~Lie
superalgebra structure on $L\otimes A$ in some cases where $L$ is a~Lie superalgebra and $A$ is a
non-supercommutative associative superalgebra.

\section{Main theorem}\label{sec:main}

\ssbegin[Main Theorem]{Theorem}[Main Theorem]\label{th} Let  $\Kee$ be an algebraically closed   field of characteristic~$p$.

\emph{1)} If $p\neq2$, any NIS on a
simple finite-dimensional Lie superalgebra is homogeneous with respect to parity, and the dimension of the space
spanned by NISes is $\leq 1$.

More precisely, the
superdimension of the space spanned by  NISes is either  $0$ (no NIS),
or $1$ \textup{(}in this case, all NISes are  even\textup{)}, or $\eps$ \textup{(}in this case, all NISes are odd\textup{)}.

\emph{2)}  If $p=2$, the superdimension of the space spanned by  NISes on a~simple
finite-dimensional Lie superalgebra is equal to either $0$ (no NIS), or $1$, or $\eps$, or $1|1$.

\textbf{This superdimension is equal
to $1|1$ if and only if the Lie superalgebra is a~queerification \emph{(see Subsection~\ref{112})} of a~simple restricted  Lie algebra with a~NIS}. \emph{(Here restrictedness is understood in the
classical sense, see Subsection~\ref{SSp2pStr};  for classification in various cases of simple Lie
(super)algebras with a~NIS, and definitions of other types of restrictedness, indigenous for $p=2$, see \cite{BKLS}.)}
\end{Theorem}

\ssec{A reformulation of Main Theorem}\label{cor} The proof of Theorem~\ref{th} in case $p=2$ actually proves the following Theorem \ref{ThP2}, which implies Theorem~\ref{th} for $p=2$ in its turn.

\sssbegin[A version of Main Theorem]{Theorem}[A version of Main Theorem]\label{ThP2} Let  $\Kee$ be an algebraically closed   field of characteristic $p=2$.

If a~simple finite-dimensional Lie  superalgebra $\fg$ has a~NIS, then

$\bullet$ the space spanned by all NIS-forms on $\fg$ is a~superspace;

$\bullet$ the dimension of the space of even IS-forms on $\fg$ is $\leq 1$;

$\bullet$ the dimension of the space of odd IS-forms on $\fg$ is $\leq 1$;

$\bullet$ i.e., the space spanned by all NIS-forms  on $\fg$ is a~superspace whose superdimension is equal to either $1$, or $\eps$, or $1|1$
\textup{(}since the even and odd components of an IS-form are also IS-forms\textup{)}.

$\bullet$ Any homogeneous IS-form on $\fg$ is either $0$ or non-degenerate.
\end{Theorem}

\begin{proof}
First, we prove several lemmas which allow us to restrict ourselves to homogenous NISes.

\ssbegin[On homogenous invariant symmetric forms for $p\neq 2$]{Lemma}[On homogenous invariant symmetric forms for $p\neq 2$]\label{homog-ne2} Let $\fg$ be a~simple Lie superalgebra over a~field of characteristics $p\neq 2$. Then, any homogenous invariant symmetric form on $\fg$ is either zero or non-degenerate.\end{Lemma}

\begin{proof} Let $\omega$ be a~homogenous IS-form on $\fg$. Then, $\Ker \omega$ is a~subsuperspace of $\fg$ invariant with respect to $\ad_\fg$, i.e., it is an ideal. Since $\fg$ is simple, either $\Ker \omega = 0$, in which case $\omega$ is non-degenerate, or $\Ker\omega = \fg$, in which case $\omega = 0$.\end{proof}

The same is not true in characteristic $2$. The problem is that even though $\Ker \omega$ is invariant with respect to $\ad_\fg$, it may be not an ideal, since it may be not closed under squaring. Before we formulate a~similar statement for $p=2$, let us prove two more lemmas, where we assume $p=2$. (The lemmas are true for $p\neq 2$ as well, but they are trivial in that case.)

\ssbegin[A necessary condition for NIS]{Lemma}[A necessary condition for NIS]\label{1} Let $p=2$. Let $\fg$ be a~simple Lie superalgebra with a~NIS $(-,-)$. Then, $[\fg,\fg]=\fg$. \end{Lemma}

\sssbegin[Peculiarity of the $p=2$ case]{Comment}[Peculiarity of the $p=2$ case]  One might think \lq\lq since $[\fg,\fg]$ is an ideal in $\fg$ which is supposed to be simple, we have nothing to prove". But for Lie \textbf{super}algebras in characteristic 2 there is a~difference between the commutant $[\fg, \fg]:=\Span([x,y]\mid x, y\in \fg)$,
and the first derived algebra. Recall that the \textit{derived Lie superalgebras} of $\fg$ are
defined to be (for $i\geq 0$)
\be\label{derive2}
\fg^{(0)}: =\fg, \quad
\fg^{(i+1)}:=\begin{cases}[\fg^{(i)},\fg^{(i)}]&\text{for $p\neq
2$,}\\
[\fg^{(i)},\fg^{(i)}]+\Span\{g^2\mid g\in\fg^{(i)}_\od\}&\text{for
$p=2$}.\end{cases}
\ee
\end{Comment}

\begin{proof}[Proof of Lemma~$\ref{1}$] Suppose $[\fg,\fg]\neq \fg$. Then, the orthogonal complement to  $[\fg,\fg]$  with respect to $(-,-)$ contains a~non-zero element~$u$. Since $\fg$ is simple, it has zero center, and hence there exists an $x\in \fg$ such that ${[u,x]\neq 0}$. Since the form $(-,-)$ is non-degenerate, there exists a~$y\in \fg$ such that $([u,x], y)\neq 0$. But then $(u, [x,y]) = ([u,x], y) \neq 0$, contradicting the fact that $u\in [\fg,\fg]^\perp$. \end{proof}

\ssbegin[A technical lemma]{Lemma}[A technical lemma]\label{2} Let $p=2$. Let $\fg$ be a~simple Lie superalgebra such that $[\fg,\fg]=\fg$, and  $S \subseteq \fg$  its subsuperspace such that $[S,\fg]\subseteq S$. Then, either $S=0$ or $S=\fg$. \end{Lemma}

\begin{proof} If $S\neq 0$, let $\overline{S}$ be the completion of $S$ with respect to squaring. Since $\overline{S}$ is a~subsuperspace of $\fg$, is closed under squaring, and  $[\overline{S}, \fg] = [S, \fg] \subseteq S$, it follows that $\overline{S}$ is an ideal, and hence $\overline{S} = \fg$. But then $\fg = [\fg, \fg] = [\overline{S}, \fg] \subseteq S$, and therefore $S=\fg$. \end{proof}

Now we can formulate the statement we will use instead of Lemma \ref{homog-ne2} when $p=2$:

\ssbegin[Analog of Lemma \ref{homog-ne2} for $p=2$]{Lemma}[Analog of Lemma \ref{homog-ne2}]\label{homog-2} Let $p=2$. Let $\fg$ be a~simple Lie superalgebra  and let there be a~NIS on $\fg$. Then, any homogenous invariant symmetric form on $\fg$ is either zero or non-degenerate.\end{Lemma}

\begin{proof} The proof is analogous to the proof of Lemma \ref{homog-ne2}, but we use Lemmas \ref{1} and \ref{2} to show that $\Ker \omega$ is either $0$ or $\fg$. \end{proof}

\textbf{Completion of the proof of Theorem~\ref{th}}.
Now we can restrict ourselves to considering homogenous non-degenerate invariant symmetric forms, since for any non-homogenous NIS, its even and odd components are invariant and symmetric, and hence non-degenerate by Lemmas \ref{homog-ne2} and \ref{homog-2}.

Fix a~basis in $\fg$, and let $B_1$ and $B_2$ be Gram matrices of \textbf{non-degenerate}  homogenous invariant symmetric forms $\omega_1$ and $\omega_2$; consider the 1-parameter family of invariant symmetric forms $\omega_\lambda$, where $\lambda\in\Kee$,  with Gram matrices $B_\lambda=B_1+\lambda B_2$. Consider $B_\lambda$ just as a~matrix, not supermatrix, and calculate its determinant; it is a~polynomial of $\lambda$.

If $\Kee$ is  \textbf{algebraically closed}, there exists a~$\lambda_0 \in\Kee$ for which ${\det B_{\lambda_0}=0}$. Then, the form $\omega_{\lambda_0}$ is degenerate. If the forms $\omega_1$ and $\omega_2$ are of the same parity, then $\omega_{\lambda_0}=0$ by Lemmas~\ref{homog-ne2} and \ref{homog-2}. This means that $\omega_1=-\lambda_0 \omega_2$, i.e., any two homogeneous  NISes of the same parity are proportional to one another. So, if $\omega_1$ and $\omega_2$ are even, then the superdimension of the space of even NISes is equal to $1$;  if $\omega_1$ and $\omega_2$ are odd, then the superdimension of the space of odd NISes is  equal to  $\eps$.

Let now forms  $\omega_1$ and $\omega_2$ be of different parity. Then, $\Ker \omega_{\lambda_0}$ is a~non-trivial $\ad_\fg$-invariant subspace, but it is not a~subsuperspace. Consider two subspaces
\[
\begin{array}{l}
V:=\Ker \omega_{\lambda_0}\cap \fg_\ev\oplus \Ker \omega_{\lambda_0}\cap \fg_\od,\\
W:=\pr_\ev(\Ker \omega_{\lambda_0})\oplus \pr_\od(\Ker \omega_{\lambda_0}),
\end{array}
\]
where $\pr_\ev$ and $\pr_\od$ are projections to $\fg_\ev$ and $\fg_\od$, respectively. Since both $V$ and $W$ are  {$\ad_\fg$-invariant} sub\textbf{super}spaces,  they are either $0$ or $\fg$ by Lemmas \ref{1} and \ref{2}. Since $\omega_{\lambda_0}$ is non-zero, $V\neq \fg$, so $V=0$; since $\omega_{\lambda_0}$ is degenerate, $W\neq 0$, so $W=\fg$.
 Hence, there exists an odd isomorphism $f$ between linear superspaces $\fg_\ev$ and $\fg_\od$ such that
\begin{equation*}\label{W}
\Ker\omega_{\lambda_0}=\Span\{x+f(x)\mid x\in\fg_\ev\}.
\end{equation*}
Since $\Ker\omega_{\lambda_0}$ is $\ad_\fg$-invariant, we see that for all $a,b\in\fg_0$:
\begin{equation}\label{fab}
\begin{array}{l}
{}[a,b+f(b)]=[a,b]+[a,f(b)]\Longrightarrow [a,f(b)]=f([a,b]); \\
{}[a+f(a),b]= [a,b]+[f(a),b]\Longrightarrow [f(a),b]=f([a,b]);
\\
{}[f(a),b+f(b)]= [f(a),b]+[f(a),f(b)]\\
{}=[f(a),f(b)]+f([a,b])\Longrightarrow f([a,b])=f([f(a),f(b)]).
\end{array}
\end{equation}

The bottom line in \eqref{fab} implies that
\be\label{this}
[f(a),f(b)]=[a,b]\text{~~for all $a,b\in \fg_\ev$.}
\ee

Up to this moment our reasoning did not depend on $p$.

Now, let $p\ne 2$.
Notice that the left-hand side of the equality~\eqref{this} is symmetric while the right-hand side is anti-symmetric. Hence,
\[
[f(a),f(b)]=[a,b]=0\text{~~and $[f(a),b]=f([a,b])=0$ for all $a,b\in \fg_\ev$},
\]
 and so $\fg$ is commutative. This contradicts the simplicity of $\fg$, and hence there can not exist two NISes of different parity on a~simple Lie superalgebra over an algebraically closed field of characteristic $p\ne 2$. There can not exist an inhomogeneous NIS in this situation either, as was mentioned above.
This completes the proof of the theorem if $p\ne 2$.

If $p=2$, no such conclusion follows  from  equality~\eqref{this}; it only tells us that $\fg$ is a~queerification (see Subsection~~\ref{112} and \cite{BLLS}) of $\fg_\ev$ and $\fg_\ev$ is (classically) restricted. Besides, the restriction of the even of the two forms~$\omega_i$ to~$\fg_\ev$ is a~NIS.

In any \textit{super}commutative superalgebra $\cA$, we have $a^2=0$ for any $a\in\cA_\od$, see condition~\eqref{supcomm}. For an example of a~commutative $\Zee/2$-graded (a.k.a. super) algebra  with a~structure of a~non-supercommutative superalgebra over $\Ree$ take $\Cee$ and set $p(i)=\od$.

Now, let  $\fg$  be a~queerification  of a~simple restricted Lie algebra $\fg_\ev$ with a~NIS $\omega$.  Then, $\fg$ can be represented in the form $\fg=\fg_\ev\otimes\cal A$, where $\cal A$ is an associative and commutative, but \textbf{not} supercommutative, 
superalgebra spanned by an even element $1$ (unit)  and an \textbf{odd} one $a$, subject to the relation $a^2=1$, and the natural bracket (for squaring, see Subsection~\ref{sqA}):
\[
{}[x\otimes \varphi, y\otimes \psi]=[x,y]\otimes \varphi\psi \text{~~for any $x,y\in\fg_\ev$, $\varphi,\psi\in\cal A$.}
\]
Determine two bilinear forms on $\fg:=\fg_\ev\otimes \cA$:
\begin{equation}\label{eq:NISonQ}
\omega_i(x\otimes \varphi, y\otimes \psi)=\omega(x,y) f_i(\varphi\psi) \text{~~for $i=1,2$},
\end{equation}
where $f_1(\alpha\cdot 1+\beta\cdot a)=\alpha$ and $f_2(\alpha\cdot 1+\beta\cdot a)=\beta$ for any $\alpha,\beta\in\Kee$.
It is clear that both these forms are non-degenerate, invariant, and symmetric; $\omega_1$ is even and $\omega_2$ is odd.
\end{proof}

\ssec{On \protect \textit{degenerate}  invariant symmetric bilinear forms}  Let $p=2$ and let $\fg$ be a~simple finite-dimensional Lie superalgebra such that $[\fg,\fg]\neq \fg$, i.e., there are elements which can not be obtained by bracketing, only by squaring, and let $k := \codim\, [\fg,\fg]$; for example, if $\fg=\fosp_{I\Pi}^{(1)}(1|2)$, then $k=2$. The dimension of the space of invariant \textbf{degenerate} symmetric bilinear forms on $\fg$ is $\geq \frac12k(k+1)=\dim S^2(\fg/ [\fg,\fg])$, because this is the dimension of the space of the  forms whose kernel contains $[\fg,\fg]$. Such forms $(-, -)$ are invariant because for them,
\[
([x, y], z) = (x, [y, z]) = 0\text{~~for any $x,y,z\in\fg$.}
\]

\section{Basics}\label{basics}

\ssec{The Sign Rule, skew and anti} The definitions of \textit{Lie superalgebra}
are the same for any $p\neq 2$ or 3: they are obtained from the definition of the Lie algebra using the Sign Rule \lq\lq if something of parity $p$ is moved past something of parity $q$, the sign $(-1)^{pq}$ accrues; formulas defined on homogeneous elements are extended to any elements via linearity".

In addition to the Sign Rule, note that
\textit{morphisms} of superalgebras are only even ones.

Observe that sometimes applying the Sign Rule requires some dexterity
(we have to distinguish between two versions both of which turn in the
nonsuper case into one, called either skew- or anti-commutativity, see \cite{Gr}):
\[
\begin{array}{ll}
ba=(-1)^{p(b)p(a)}ab &(\text{super commutativity})\\
ba=-(-1)^{p(b)p(a)}ab &(\text{super anti-commutativity})\\
ba=(-1)^{(p(b)+1)(p(a)+1)}ab &(\text{super skew-commutativity})\\
ba=-(-1)^{(p(b)+1)(p(a)+1)}ab
&(\text{super antiskew-commutativity})\end{array}
\]
 In the case of characteristic 2, of main interest to us, the above conditions turn into one, so  \textit{super commutativity} and  \textit{super antiskew-commutativity} are defined as
\be\label{supcomm}
\text{$ab=ba$ for all $a,b$, and $a^2=0$ for $p(a)=\od$,}
\ee
whereas  \textit{super anti-commutativity} and  \textit{super skew-commutativity}  are defined as
\[
\text{$ab=ba$ for all $a,b$, and $a^2=0$ for $p(a)=\ev$.}
\]

\ssec{Lie superalgebra, pre-Lie superalgebra, Leibniz superalgebra}\label{sring} Generalization of notions of this Subsection to  super rings and modules over them is immediate.

\underline{For any $p$}, a~\textit{Lie superalgebra} is a
superspace $\fg=\fg_\ev\oplus\fg_\od$ such that the even part
$\fg_\ev$ is a~Lie algebra, the odd part $\fg_\od$ is a
$\fg_\ev$-module (made into the two-sided one by
\textit{anti}-symmetry, i.e., $[y,x]=-[x,y]$ for any $x\in
\fg_\ev$ and $y\in
\fg_\od$)  and on $\fg_\od$, a~\textit{squaring} $x\mapsto x^2$ and the \textit{bracket} are defined via a~linear map $s: S^2(\fg_\od)\tto\fg_\ev$, where  $S^2$ denotes the operator of raising to symmetric square, as follows, for any $x,y\in\fg_\od$:
\begin{eqnarray}\label{squaring0}
x^2 := s(x\otimes x);\label{squaring1}\\
{}[x,y] := s(x\otimes y + y\otimes x).
\label{squaring2}
\end{eqnarray}

The linearity of the $\fg_\ev$-valued function $s$ implies that
\begin{eqnarray}\label{squaring3}
\text{$(ax)^2=a^2x^2$ for any $x\in
\fg_\od$ and $a\in \Kee$, and}\label{squaring4}\\
{}[\cdot,\cdot]\text{~is a~bilinear form on $\fg_\od$ with values
in $\fg_\ev$.}\label{squaring5}
\end{eqnarray}

The \textit{Jacobi identity} involving odd elements takes the form of the following two conditions:
\begin{eqnarray}
~[x^2,y]=[x,[x,y]]\text{~for any~} x\in\fg_\od, y\in\fg_\ev,\label{6b}\\
~[x^2,x]=0\;\text{ for any $x\in\fg_\od$.}\label{6}
\end{eqnarray}

In the literature describing Lie superalgebras in characteristics $p\neq 2,3$, the Jacobi identity is usually given for three elements, and when all those elements are odd, it takes the form
\begin{equation}\label{Jac_3odd}
~[x,[y,z]]+[y,[z,x]]+[z,[x,y]]=0\text{ for any $x,y,z\in\fg_\od$.}
\end{equation}
This equation follows from \eqref{6} in any characteristic; also, if $p\neq 2,3$, then \eqref{6} follows from \eqref{Jac_3odd}, since, by substituting $x=y=z$, we get $6[x,x^2]=0$. The more general equation \eqref{6} is necessary for the super version of the Poincare-Birkhoff-Witt theorem to hold, and this is why we use it.

\underline{If $p=3$}, and  the Jacobi identity is formulated in the usual way, i.e., for three elements, then
the condition
\be\label{p=3}
[x,[x,x]]=0 \text{~~for any $x\in\fg_\od$}
\ee
should be added to the antisymmetry and Jacobi identity amended by the Sign Rule, separately.  If $p=3$, the superalgebra satisfying the antisymmetry and the Jacobi identity, but not the condition~\eqref{p=3} is called \textit{pre-Lie superalgebra}; for examples, see \cite{BeBou}.

\underline{If $p=2$}, the \textit{antisymmetry} for $p=2$ should be replaced by an equivalent for $p\neq 2$, but otherwise
stronger \textit{alternating} or \textit{antisymmetry}    condition
\[
[x,x]=0\text{~~for any $x\in\fg_\ev$}.
\]

For any $p$, the superalgebra satisfying the Jacobi identity, and without any restriction on symmetry is called a~\textit{Leibniz superalgebra}.

Over $\Zee/2$, the condition \eqref{6}  must (see Example \ref{exa}) be replaced with a~more general one:
\be\label{JIgen}
\text{$[x^2,y]=[x,[x,y]]\text{~for any~} x,y\in\fg_\od$.}
\ee
For any other ground field this condition is equivalent to condition \eqref{6}.

\sssbegin[On the Jacobi identity over $\Zee/2$]{Example}[On the Jacobi identity over $\Zee/2$]\label{exa} This example shows that over $\Zee/2$, condition \eqref{JIgen} is not a~corollary of condition \eqref{6}. Take a~$2|3$-dimensional superalgebra $\fg$ with the even part spanned by elements $A$ and $B$, the odd part spanned by elements $X$, $Y$ and $Z$, and the algebraic structure given as follows:

$\bullet$ the even part is commutative;

$\bullet$ the action of $\fg_\ev$ on $\fg_\od$ is given by the following multiplication table
\[
\begin{tabular}{|l |l |l |l |}
\hline
&$X$&$Y$&$Z$\\
\hline
$A$&0&$Z$&0\\
\hline
$B$&$Z$&0&0\\
\hline
\end{tabular}
\]
$\bullet$ the squaring on $\fg_\od$ is given by the formula
\[
(aX + bY + cZ)^2 = a^2A + b^2B\text{~for all~} a,b,c\in \Zee/2.
\]

This algebra $\fg$ satisfies (\ref{6}) since
\[
{}[(aX + bY + cZ)^2, aX + bY + cZ] = [a^2A + b^2B, aX + bY + cZ] = (a^2b+ab^2)Z,
\]
and since $a,b\in \Zee/2$, we have $a^2=a$ and $b^2=b$, so $a^2b+ab^2=0$.

But the algebra $\fg$ does not satisfy condition (\ref{JIgen}): 
\[
{[X^2, Y] = Z \neq [X, [X, Y]] = 0}.
\]
\end{Example}

\sssec{Ideals, simplicity, derived algebras, modules}
By an \textit{ideal} of a~Lie superalgebra one always means an \textit{homogeneous} ideal; for $p=2$, the ideal should be closed with respect to squaring.

The Lie superalgebra $\fg$ is said to be \textit{simple} if $\dim \fg >1$ and $\fg$ has no nontrivial (distinct from 0 and $\fg$) ideals.

For $p=2$, the definition of the derived algebra of the Lie superalgebra $\fg$ changes, see eq.~\eqref{derive2}.

An even linear map $r: \fg\tto\fgl(V)$ is said to be a
\textit{representation of the Lie superalgebra} $\fg$ and $V$ is a
\textit{$\fg$-module} if
\begin{equation*}\label{repres}
\begin{array}{l}
r([x, y])=[r(x), r(y)]\quad \text{ for any $x, y\in
\fg$;}\\
r(x^2)=(r(x))^2\text{~for any $x\in\fg_\od$.}
\end{array}
\end{equation*}

Since we want $\fder\ \fg$ to be a~Lie superalgebra for any Lie superalgebra $\fg$, we have to add a~generalization of conditions \eqref{6b}, \eqref{6}, namely, the condition
\be\label{9.10}
D(x^{2}) = [D(x),x]\text{~~for odd elements $x\in\fg$ and any $D\in\fder\ \fg$}.
\ee
Clearly, condition \eqref{9.10} turns into conditions \eqref{6b}, \eqref{6} for $D=\ad_y$, where $y\in\fg$.

\subsection{The $p|2p$-structure or restricted Lie
superalgebra}\label{SSp2pStr} Let the ground field~$\Kee$
be of characteristic $p>0$, let $\fg$ be a~Lie algebra. For every $x\in
\fg$, the operator $(\ad_x)^{p}$ is a~derivation of
$\fg$. If this derivation is an inner one, then there is a~map
(called \textit{$p$-structure}) ${}[p]:\fg\tto\fg, \ x\mapsto x^{[p]}$ such that
\begin{equation*} \label{restricted-3}
\begin{array}{l}
{}[x^{[p]}, y]=(\ad_x)^{p}(y)\quad \text{~for any~}x,y\in\fg,\\
(ax)^{[p]}=a^px^{[p]}\quad \text{~for any~}a\in\Kee,~x\in\fg,\\
(x+y)^{[p]}=x^{[p]}+y^{[p]}+\mathop{\sum}\limits_{1\leq i\leq
p-1}s_i(x, y) \quad\text{~for any~}x,y\in\fg,\\
\end{array}
\end{equation*}
where $is_i(x, y)$ is the coefficient of $\lambda^{i-1}$ in
$(\ad_{\lambda x+y})^{p-1}(x)$,
then the Lie algebra
$\fg$ is said to be \emph{restricted} or \emph{having a
$p$-structure}.

For a~Lie superalgebra $\fg$ in characteristic $p>0$, let the Lie
algebra $\fg_\ev$ be restricted and
\begin{equation} \label{restr3}
[x^{[p]}, y]=(\ad_x)^{p}(y)\quad \text{~for any~}x\in\fg_\ev,~y\in\fg.
\end{equation}
This gives rise to the map
\begin{equation*}\label{2p}
{}[2p]:\fg_\od\to\fg_\ev, ~~~x\mapsto(x^2)^{[p]},
\end{equation*}
satisfying the condition
\begin{equation}\label{2p2}
{}[x^{[2p]},y]=(\ad_x)^{2p}(y)\quad\text{~for any~}x\in\fg_\od,~y\in\fg.
\end{equation}
The pair of maps $[p]$ and $[2p]$ is called a~$p$-\textit{structure}
(or, sometimes, a~$p|2p$-\textit{structure}) on~$\fg$, and~$\fg$ is
said to be \textit{restricted}, cf. \cite{WZ}. It suffices to determine the $p|2p$-structure on any basis of $\fg$; on simple Lie superalgebras there is at most one $p|2p$-structure.

$\bullet$ If condition \eqref{restr3} is not satisfied,
the $p$-structure on $\fg_\ev$ does not have to generate a
$p|2p$-structure on~$\fg$: even if the actions of $(\ad_x)^p$ and
$\ad_{x^{[p]}}$ coincide on $\fg_\ev$, they do not have to coincide
on the whole of $\fg$. 



$\bullet$ If $p=2$, there are \textbf{other}, no less natural, versions of restrictedness, see \cite{BLLS}; we will not consider them in this text.

\subsection{Linear (matrix) Lie (super)algebras}\label{SS:2.3}
The \textit{general linear} Lie
superalgebra of all linear operators in the superspace $V=V_{\bar 0}\oplus V_{\bar 1}$ over the ground field $\Kee$ is denoted by
$\fgl(\Size)$, where  $\Size=(p_1, \dots, p_{|\Size|})$ is an ordered
collection of parities of the basis vectors of $V$ for which we take only vectors \textit{homogeneous with respect to parity} and $|\Size|:=\dim V$. The  Lie
superalgebra of all supermatrices of size $\Size=(p_1, \dots, p_{|\Size|})$ is also denoted by
$\fgl(\Size)$.
Usually, for the \textit{standard} (simplest from a~certain point of view) format $\Size_{st}:=(\ev,
\dots, \ev, \od, \dots, \od)$, the notation $\fgl(\Size_{st})$ is abbreviated to $\fgl(\dim V_{\bar
0}|\dim V_{\bar 1})$. Any $X\in\fgl(\Size)$ can
be uniquely expressed as the sum of its even and odd parts; in the
standard format this is the following block expression; on non-zero summands the parity is defined:
\[
X=\mmat A,B,C,D,=\mmat A,0,0,D,+\mmat 0,B,C,0,,\quad
 p\left(\mmat A,0,0,D,\right)=\ev, \; p\left(\mmat 0,B,C,0,\right)=\od.
\]

The \textit{supertrace} is the map $\fgl (\Size)\tto \Kee$,
$(X_{ij})\longmapsto \sum (-1)^{p_{i}(p(X)+1)}X_{ii}$.

Thus, in the standard format, $\str \begin{pmatrix}A&B\\ C&D\end{pmatrix}=\tr A- \tr D$.
Observe that for Lie superalgebra $\fgl_\cC(p|q)$ over a~supercommutative superalgebra $\cC$, i.e., for supermatrices with elements in $\cC$, we have
\be\label{pX}
\begin{array}{l}
\str X=\tr A- (-1)^{p(X)}\tr D\text{~~for any $X=\begin{pmatrix}A&B\\ C&D\end{pmatrix}$,}\\
\text{where $p(X)=p(A_{ij})= p(D_{kl})=p(B_{il})+\od=p(C_{kj})+\od$},\\
\end{array}
\ee
so  if $\cC_\od\neq 0$, the supertrace coincides with the trace on odd supermatrices.

Since $\str\ [x, y]=0$, the subsuperspace of supertraceless
matrices constitutes a~Lie subsuperalgebra called \textit{special linear} and denoted
$\fsl(\Size)$.

There are, however, at least two super versions of $\fgl(n)$, not
one; for reasons, see \cite[Ch1, Ch.7]{Lsos}. The other version --- $\fq(n)$ --- is called the \textit{queer}
Lie superalgebra and is defined as the one that preserves --- if $p\neq 2$ --- the
complex structure given by an \textit{odd} operator $J$, i.e.,
$\fq(n)$ is the centralizer $C(J)$ of $J$:
\[
\fq(n)=C(J)=\{X\in\fgl(n|n)\mid [X, J]=0 \}, \text{ where }
J^2=-\id.
\]
It is clear that by a~change of basis we can reduce $J$ to the form (shape)
$J_{2n}:=\mat {0&1_n\\-1_n&0}$
in the standard
format, and then $\fq(n)$ takes the form
\begin{equation}\label{q}
\fq(n)=\left \{(A,B):=\mat {A&B\\B&A}, \text{~~where $A, B\in\fgl(n)$}\right\}.
\end{equation}
(Over any algebraically closed field $\Kee$, instead of $J$ we can take any odd operator $K$ such that $K^2=a\id_{n|n}$, where $a\in \Kee^\times$; and the centralizers of $K$, Lie superalgebras $C(K)$, are isomorphic for distinct $K$; if $p=2$, it is natural to select $K^2=\id$.)

The supertrace vanishes on $\fq(n)$, but there is defined the odd \textit{queertrace}  $\qtr\colon (A,B)\longmapsto \tr B$ which vanishes on the first derived of $\fq(n)$, so it \textbf{is} a~trace. Denote by $\fsq(n)$ the Lie superalgebra
of \textit{queertraceless} matrices; set $\fp\fsq(n):=\fsq(n)/\Kee 1_{2n}$.

If $p=2$, there is a~tailor-made trace on $\fsq(n)$; it is not a~trace on $\fq(n)$. 
It is the \textit{halftrace} ($\frac12\str$) given by the formula 
\[
\htr (A, B):=\tr A;
\]
we set
\[\begin{array}{l}
\fs_e\fsq(n):=\{X=(A,B)\in \fsq(n)\mid \tr(A)=0\};\\
\fp\fs_e\fsq(2n):=\fs_e\fsq(2n)/\Kee 1_{4n}.
\end{array}
\]

Clearly, $\fgl$ and $\fq$ correspond to the super version of Schur's lemma over an algebraically closed field: an irreducible module over a~collection $S$ of homogeneous operators can be \textit{absolutely irreducible}, i.e., have no proper invariant subspaces, in which case the only operator commuting with $S$ is a~scalar (the $\fgl$ case), or can have an invariant subspace which is not a~sub\textbf{super}space, in which case the superdimension of the module is of the form $n|n$ and an odd operator $J$ interchanges the homogeneous components of the module (the $\fq$ case).

\subsubsection{Supermatrices  of bilinear forms}\label{sssBilMatr}
Following~\cite[\S1.5]{L2}, \cite[Ch.3, \S5]{Ma}, and~\cite[Ch.1]{Lsos},
to any linear mapping $F: V\tto W$ of superspaces there corresponds the
dual mapping $F^*:W^*\tto V^*$ between the dual superspaces. In a~basis
consisting of the vectors $v_{i}$ of format $\Size$, the formula
$F(v_{j})=\sum_{i}v_{i}X_{ij}$ assigns to $F$ the
supermatrix $X$. In the dual bases, the \textit{supertransposed}\index{Supertransposition}
matrix $X^{st}$ corresponds to $F^*$:
\[
(X^{st})_{ij}=(-1)^{(p_{i}+p_{j})(p_{i}+p(A))}X_{ji}.
\]
In the standard format, this means
\[
X=\mmat A,B,C,D,\longmapsto\begin{cases}\begin{pmatrix}A^t&C^t\\-B^t&D^t\end{pmatrix}&\text{if $p(X)=\ev$},\\
\begin{pmatrix}A^t&-C^t\\B^t&D^t\end{pmatrix}&\text{if $p(X)=\od$}.\\
\end{cases}
\]

The supermatrices $X\in\fgl(\Size)$ such that
\[
X^{st}B+(-1)^{p(X)p(B)}BX=0\quad \text{for an homogeneous matrix
$B\in\fgl(\Size)$}
\]
constitute the Lie superalgebra $\faut (B)$ that preserves the
bilinear form $\cB$ on $V$ whose Gram matrix $B=(B_{ij})$ is given by the formula 
\be\label{martBil}
B_{ij}=(-1)^{p(B)p(v_i)}\cB(v_{i}, v_{j})\text{~~for the basis vectors $v_{i}\in V$}
\ee
in order to identify a~bilinear form $B(V, W)$ with an element of $\Hom(V, W^*)$. 
Consider the \textit{upsetting} of bilinear forms
$u\colon\Bil (V, W)\tto\Bil(W, V)$ given by the formula \be\label{susyB}
u(\cB)(w, v)=(-1)^{p(v)p(w)}\cB(v,w)\text{~~for any $v \in V$ and $w\in W$.}
\ee
In terms of the Gram matrix $B$ of $\cB$: the form
$\cB$ is  \textit{symmetric} if  and only if 
\be\label{BilSy}
u(B)=B,\;\text{ where $u(B)=
\mmat R^{t},(-1)^{p(B)}T^{t},(-1)^{p(B)}S^{t},-U^{t},$ for $B=\mmat R,S,T,U,$.}
\ee
Similarly, \textit{anti-symmetry} of  
$\cB$ means that $u(B)=-B$. There are no ``supersymmetric" bilinear forms; in \cite{BKLS} this is a~typo, NISes were meant.


Observe
that \textbf{the passage from $V$ to $\Pi (V)$ turns every symmetric
form $\cB$ on $V$ into an anti-symmetric one $\cB^\Pi$ on $\Pi (V)$  and anti-symmetric $\cB$ into symmetric $\cB^\Pi$  by setting}
\[
\text{$\cB^\Pi(\Pi(x), \Pi(y)):=(-1)^{p(B)+p(x)+p(x)p(y)}\cB(x,y)$ for any $x,y\in V$}.
\]

If $p\neq 2$, the Lie superalgebra preserving a~non-degenerate symmetric even bilinear form $\cB$ on the $a|2b$-dimensional superspace  is denoted $\fosp_\cB(a|2b)$ or just $\fosp(a|2b)$ (for a basis in a standard format); if $p= 2$, it is denoted $\fo\fo_\cB(a|b)$, see \cite{BGL1}.

\subsection{Queerification for $p=2$ (from \cite{BLLS})}\label{112} If $p=2$, then we
can queerify any restricted Lie algebra $\fg$ as follows. We set
$\fq(\fg)_\ev=\fg$ and $\fq(\fg)_\od=\Pi(\fg)$; define the
multiplication involving the odd elements as follows:
\begin{equation}\label{q(g)}
{}[x,\Pi(y)]=\Pi([x,y]);\quad (\Pi(x))^2=x^{[2]}\quad\text{~for
any~~}x,y\in\fg.
\end{equation}
Clearly, if $\fg$ is restricted and $\fii\subset\fq(\fg)$ is an ideal,
then $\fii_\ev$ and $\Pi(\fii_\od)$ are ideals in $\fg$. So, if $\fg$ is
restricted and simple, then $\fq(\fg)$ is a~simple Lie superalgebra.
(Note that $\fg$ has to be simple as a~Lie algebra, not just as a
\textit{restricted} Lie algebra, i.e., $\fg$ is not allowed to have
\textbf{any} ideals, not only restricted ones.) A generalization of
the queerification is the following procedure producing as many
simple Lie superalgebras as there are simple Lie algebras.

\subsubsection{Generalized queerification}\label{Method1} Let the \textit{$1$-step
restricted closure} $\fg^{<1>}$ of the  Lie algebra~$\fg$ be
the minimal subalgebra of the (classically) restricted closure
$\overline{\fg}$ containing $\fg$ and all the elements $x^{[2]}$,
where $x\in\fg$. To any
Lie algebra $\fg$ the \textit{generalized queerification} assigns the Lie superalgebra
\begin{equation*}\label{tildeQ}
\tilde \fq(\fg):=\fg^{<1>}\oplus \Pi(\fg)
\end{equation*}
with squaring given by $(\Pi(x))^2=x^{[2]}$ for any $x\in\fg$.
Obviously, for $\fg$ restricted, the generalized queerification coincides with the
queerification: $\tilde \fq(\fg)=\fq(\fg)$. As proved in  \cite{BLLS}, if $\fg$ is a~simple Lie algebra, then $\tilde \fq(\fg)$ is a~simple Lie superalgebra.

\subsection{Vectorial Lie superalgebras }\label{ssVectorial} For their definition and NISes on them, see \cite{BKLS}.

\ssec{Deformations with odd parameters and deforms} Which of the infinitesimal deformations can be extended to
a~global one is a~separate much tougher question, usually solved
\textit{ad hoc}; for examples over fields of characteristics $3$ and
$2$, see \cite{BLW} and references therein. Deformations with odd parameters are always integrable.
Let us give two graphic examples.

1) \textbf{Deformations of representations}.
Consider a~representation $\rho:\fg\tto\fgl(V)$.
The tangent space of the moduli superspace of deformations of $\rho$
is isomorphic to $H^1(\fg; V\otimes V^*)$. For example, if $\fg$ is
the $0|n$-dimensional (i.e., purely odd) Lie superalgebra (with the
only bracket possible: identically equal to zero), its only
irreducible representations are the 1-dimensional trivial one,
$\mathbbmss{1}$, and $\Pi(\mathbbmss{1})$. Clearly,
\[
\mathbbmss{1}\otimes \mathbbmss{1}^*\simeq
\Pi(\mathbbmss{1})\otimes \Pi(\mathbbmss{1})^*\simeq \mathbbmss{1},
\]
and, because the Lie superalgebra $\fg$ is commutative, the
differential in the cochain complex is zero. Therefore
\[
H^1(\fg; \mathbbmss{1})=E^1(\fg^*)\simeq\Pi(\fg^*),
\]
so there are $\dim\ \fg$ odd parameters of deformations of the
trivial representation.

2) \textbf{Deformations of the brackets}.
Let $\cC$ be a~finitely
generated supercommutative superalgebra, let $\Spec \cC$ be the affine super scheme defined literally as the affine scheme of any commutative ring, see \cite{L1}.

A~\textit{deformation} of
a Lie superalgebra $\fg$ over $\Spec \cC$
is a~Lie superalgebra $\fG$ over $\cC$ such that for some closed point $\fp_I\in \Spec \cC$ corresponding to  a~maximal ideal in $I\subset \cC$, we have $\fG\otimes_I\Kee\simeq\fg$ as Lie superalgebras over $\Kee$. Note that since $\cC/I\simeq\Kee$ as $\Kee$-algebras,
this definition implies that $\fG\simeq\fg\otimes \cC$ as $\cC$-modules. The deformation is \textit{trivial} if $\fG\simeq\fg\otimes \cC$ as Lie superalgebras over $\cC$, not just as modules, and \textit{non-trivial} otherwise.  (We superized \cite{Ru}, where the non-super case was considered.)

Generally, the \textit{deforms} --- the results of a~deformation --- of
a~Lie superalgebra $\fg$ over $\Kee$ are Lie superalgebras $\fG\otimes_{I'}\Kee$, where $\fp_{I'}$ is any closed point in $\Spec \cC$. 

In particular, consider a~deformation with an odd parameter $\tau$ of a~Lie superalgebra $\fg$ over field $\Kee$. This is a~Lie superalgebra $\fG$ over $\Kee[\tau]$ such that $\fG\otimes_I\Kee\simeq\fg$, where $I$
is the only maximal ideal of $\Kee[\tau]$. This implies that $\fG$ isomorphic to $\fg\otimes\Kee[\tau]$ as a~\textbf{module over $\Kee[\tau]$};  if, moreover, $\fG=\fg\otimes\Kee[\tau]$ as a~\textbf{Lie superalgebra over $\Kee[\tau]$}, i.e.,
\[
[a\otimes f, b\otimes g]=(-1)^{p(f)p(b)}[a,b]\otimes fg \text{~~for all $a,b\in \fg$ and $f,g\in\Kee[\tau]$},
\]
then the deformation is considered \textit{trivial} (and \textit{non-trivial} otherwise). 
Observe that $\fg\otimes \tau$ is not an ideal of $\fG$: any ideal should be a~free $\Kee[\tau]$-module.

\sssec{Comment: if deformations are with odd  parameters or even but formal, there are no ``deforms"} In a~sense, the people who ignore odd parameters of deformations have a~point: we (rather they) consider classification of simple Lie superalgebras (or whatever other problem) over the ground field $\Kee$, right? Not quite. Actually, the odd parameters of deformations are no less natural than the odd part of the  Lie superalgebra itself. However, to see these parameters, we have to consider whatever we are deforming not over $\Kee$, but over $\Kee[\tau]$.

We do the same when $\tau$ is even and we consider formal deformations over $\Kee[[\tau]]$. If
the formal series in $\tau$ converges in a~domain $D$, we can evaluate $\tau$ for any $\tau\in D$ and consider copies $\fg_\tau$ of $\fg$, where $\tau\in D$.  
If the parameter is formal or odd, such an evaluation is possible only trivially: $\tau\mapsto 0$.

\section{Examples}\label{sec:examples}

For descriptions of Lie superalgebras with indecomposable Cartan matrix in characteristic $p>0$ some of which can be deformed and used below, see \cite{BGL1}. For examples (even classification in several cases) of deforms of known symmetric simple modular Lie superalgebras, see \cite{BGL2}, where the cocycles we consider below are given explicitly, in terms of a~Chevalley basis. Here we consider one of the simplest examples of deforms with an odd parameter and several other examples.

\ssbegin[No NIS on $\fo\fo_{I\Pi}^{(1)}(1|2)$ and $\fo\fo_{II}^{(1)}(1|2)$]{Lemma}[No NIS on $\fo\fo_{I\Pi}^{(1)}(1|2)$ and $\fo\fo_{II}^{(1)}(1|2)$] Consider the Lie superalgebra $\fo\fo_{I\Pi}^{(1)}(1|2)$ and its deform with the help of the cocycle $c_{-2}$, and $\fo\fo_{II}^{(1)}(1|2)$ and its deforms with the help of the cocycles $c_1$ or $c_2$, see~\cite[\S 7.2
and \S 7.3]{BGL2}. There is no NIS
on any of these deformed\  Lie superalgebras $\fg:=\fg_{c_i}$.
\end{Lemma}
\begin{proof}
The direct computations show that $[\fg,\fg]\neq\fg$. Hence, no NIS on $\fg$ due to Lemma~\ref{1}.
\end{proof}

\ssec{On a~map sending cochains of~$\textbf{F}(\fg)$ to cochains of~$\fg$}\label{ssec:cochainsMap} If $p=2$,
and $\fg$ is a~Lie superalgebra, let~$\textbf{F}(\fg)$ be the desuperization of~$\fg$, i.e., \textbf{F} is
the functor that forgets squaring and parity. The forgetful  functor gives a~$\Kee$-linear map $i\colon \fg\tto\textbf{F}(\fg)$, and with its help the choice of a~basis in~$\fg$ induces the choice of a~basis in~$\textbf{F}(\fg)$.
Actually, $i$ is just the identity map on $\fg$ as a~vector space, and a~basis of $\textbf{F}(\fg)$ consists of the same vectors as a~basis in $\fg$.

If $p=2$, then $E^{\bcdot}(V)\subsetneq S^{\bcdot}(V)$, where $E^{\bcdot}(V)$ is the exterior algebra and $S^{\bcdot}(V)$ is the symmetric algebra of the space $V$. The map~$i$ induces an injective map
${i_\ast\colon C^{\bcdot}(\textbf{F}(\fg))\tto C^{\bcdot}(\fg)}$ between spaces of cochains.
The map~$i_\ast$ does not necessary commute with the differential: as it was noted in~\cite[\S\,7.11.1]{BGL2}, not every cocycle
of~$\textbf{F}(\fg)$ defines a~cocycle of~$\fg$.

Interestingly, the map~$i_\ast$ sometimes allows us to express some of cocycles of $\fg$,
representing cohomology classes, in terms of the images of cocycles of~$\textbf{F}(\fg)$ under $i_\ast$ (plus,
perhaps, terms defining a~deformation of the squaring), see examples in Lemmas \ref{Le1} and \ref{Le2}.

\ssec{Two lemmas about NISes on deforms}\label{2lemmas} The NISes on deforms of~$\fwk(n;\alpha)$, where $n=3$ or 4,  can be directly translated to the corresponding superizations --- deforms
of~$\fbgl(n;\alpha)$ --- because the squaring is not involved at all in the invariance condition for the bilinear form.

\sssbegin[On NISes on $\fbgl(4;\alpha)$]{Lemma}[On NISes on $\fbgl(4;\alpha)$]\label{Le1}
For Lie superalgebra $\fg=\fbgl(4;\alpha)$, where $\alpha\neq0,1$, all deforms depend on even parameters, see~\cite{BGL2}.

Choose a~basis in~$\fg$ as in \cite{BGL2}.
These deforms of~$\fg$ preserve a~NIS with the same Gram matrix as that of the NIS on~$\fg$, except for the deform~$\fg_{c_0}$
of~$\fbgl(4;\alpha)$ with cocycle $c_0$ when the Gram matrix~$\Gamma_{c_0}$ is a~different one. The
desuperization of~$\Gamma_{c_0}$ coincides with the Gram matrix, described in~\cite[Claim 3.3]{BKLS}, of the
corresponding deform of~$\fwk(4;\alpha)$ with cocycle~$c_0$.
\end{Lemma}

\begin{proof}
For a~basis in~$H^2(\textbf{F}(\fg);\textbf{F}(\fg))$ take (the classes of) $c_{\pm12}$, $c_{\pm10}$,
$c_{\pm8}^1$, $c_{\pm8}^2$, $c_{\pm6}^1$, $c_{\pm6}^2$, $c_{\pm4}^1$, $c_{\pm4}^2$, $c_{\pm2}$, $c_{0}$; for
their explicit expressions, see~\cite[\S~7.12]{BGL2}. Due to the symmetry of the root system, it suffices to
consider only cocycles of non-positive degree. As it was noted in~\cite[\S~7.13]{BGL2}, the corresponding basis (of non-positive degree) for~$H^2(\fg;\fg)$ consists of the following cocycles
 \footnotesize{\begin{align*}
 \bar c_{-12} ={}& i_\ast(c_{-12}),&\bar c_{-10} ={}& i_\ast(c_{-10}), \\
 \bar c_{-8}^{1} ={}& i_\ast(c_{-8}^{1}),&
 \bar c_{-8}^{2} ={}& i_\ast(c_{-8}^{2}) + \alpha^2(1+\alpha) h_3 \otimes (\hat x_{11})^{\wedge 2},\\
 \bar c_{-6}^{1} ={}& i_\ast(c_{-6}^{1}) + (1+\alpha) h_3 \otimes (\hat x_9)^{\wedge 2},&
 \bar c_{-6}^{2} ={}& i_\ast(c_{-6}^{2}) + \alpha^2(1+\alpha)(h_3 + h_4)\otimes (\hat x_8)^{\wedge 2},\\
 \bar c_{-4}^{1} ={}& i_\ast(c_{-4}^{1}) + (1+\alpha)(h_3 + h_4)\otimes (\hat x_6)^{\wedge 2},&
 \bar c_{-4}^{2} ={}& i_\ast(c_{-4}^{2}) + (1+\alpha)(h_3 + h_4)\otimes (\hat x_6)^{\wedge 2},\\
 \bar c_{-2}={}& i_\ast(c_{-2}) + (1+\alpha) h_4 \otimes (\hat x_1)^{\wedge 2},&
 \bar c_0={}& i_\ast(c_0).
 \end{align*}
}\normalsize
where $x_i$, $h_i$, $y_i$ are elements of the Chevalley basis of~$\fg$, see~\cite{BGL1}, and $\hat x_i$,
$\hat{h}_i$, $\hat y_i$ are the elements of the corresponding dual basis of~$\Pi(\fg^*)$.
\end{proof}

\sssbegin[On NISes on $\fbgl^{(1)}(3;\alpha)/\fc$]{Lemma}[On NISes on $\fbgl^{(1)}(3;\alpha)/\fc$]\label{Le2}
All deforms of the Lie superalgebra $\fg=\fbgl^{(1)}(3;\alpha)/\fc$, where $\alpha\neq0,1$,  depend on even parameters;
see~\cite{BGL2}. Choose a~basis in~$\fg$  as in \cite{BGL2}.

These deforms preserve a~NIS with the same Gram matrix as that of the NIS on~$\fg$, except for the deform~$\fg_{c_0}$
of~$\fbgl^{(1)}(3;\alpha)/\fc$ with cocycle $c_0$ when the Gram matrix~$\Gamma_{c_0}$ is a~different one. The
desuperization of~$\Gamma_{c_0}$ coincides with the Gram matrix, described in~\cite[Claim 3.4]{BKLS}, of the
corresponding deform of~$\fwk^{(1)}(3;\alpha)/\fc$ with cocycle~$c_0$.
\end{Lemma}

\begin{proof}
For the basis in~$H^2(\textbf{F}(\fg);\textbf{F}(\fg))$ take (the classes of) $c_{\pm6}$,
$c^1_{\pm4}$, $c^2_{\pm4}$, $c_{\pm2}$, $c_{0}$; for their explicit expressions;
see~\cite[\S~7.10]{BGL2}. Due to symmetry of the root system, it suffices to consider only cocycles of
non-positive degree. The corresponding basis (of non-positive degree) in $H^2(\fg;\fg)$ consists of the following cocycles, see~\cite[\S~7.11]{BGL2}
  \[
    \bar c_{-4} = i_\ast(c^2_4),\quad \bar c_{-2} = i_\ast(c_{-2}), \quad \bar c_{0} = i_\ast(c_{0}).\hfill\qed
  \]
\noqed\end{proof}

\ssbegin[NISes on $\fq(\fsl(3))$]{Claim}[NISes on $\fq(\fsl(3))$]\label{NISqsl3} For $\fg=\fq(\fsl(3))$,
consider its two deformations: $\fg_{A}(\tau)$ given by the odd cocycle~$A$ with odd parameter $\tau$ and $\fg_B(\xi)$ with parameter $\xi\in\Kee$ given by the even
cocycle~$B$;  for explicit expressions of these cocycles, see~\cite[Lemma 8.3]{BGL2}.

The space of NISes on~$\fg_A(\tau)$ is of rank $1|1$
over~$\Kee[\tau]$.

The spaces of NISes on~$\fg_B(\xi)$ is $(1|1)$-dimensional if $\xi \neq 1$, and $0$-dimensional
if $\xi = 1$.
\end{Claim}

\begin{proof} Although the proof is analytical, the above is called Claim because the expressions of NISes are obtained with the aid of \textit{SuperLie}.

Choose a~basis in~$\fg$ given by the
Chevalley basis in~$\fsl(3)$, namely, $x_i$, $y_i$, $h_i$, and $\Pi x_i$, $\Pi y_i$, $\Pi h_i$.
By Theorem~\ref{th}, there are two NISes, $\omega_{\ev}$ and $\omega_{\od}$, on~$\fq(\fsl(3))$ induced by the
even NIS on~$\fsl(3)$, see \cite{BKLS}, given by the formula~\eqref{eq:NISonQ}.

A) For the odd cocyle~$A$, a~NIS~$\omega$ on~$\fg_{A}(\tau)$ is defined as follows
\[
   \omega = a_1 \omega_{\ev} + a_2 \omega_{\od}
  + \tau (a_3 \omega_{\ev} + a_4 \omega_{\od} + a_1 B_1 + a_2 B_2),
  \text{~where~} a_1,a_2,a_3,a_4 \in \Kee,
\]
where $B_1$ (resp.  $B_2$) is an odd (resp. even) bilinear form for which
\[
\begin{array}{l}
  B_1(\Pi h_1,h_1) =   B_1(\Pi h_1,h_2) =   B_1(\Pi h_2,h_1) =  B_1(\Pi h_2,h_2) = 1, \\
  B_2(y_1,x_1) = B_2(y_3,x_3) = B_2(\Pi h_2,\Pi h_1) = B_2(\Pi y_2,\Pi x_2) = 1.
  \end{array}
\]
and zero on all other pairs of Chevalley basis.  Observe that~$\fg_A(\tau)$ is a~free module over~$\Kee[\tau]$. An
arbitrary element $t = a~+ b\tau\in\Kee[\tau]$ is defined by a~pair of numbers $a,b\in\Kee$. Set
$t_1=a_1 + a_3\tau$ and $t_2 = a_2 + a_4\tau\in \Lambda[\tau]$. We have
\[
  \omega = t_1(\omega_{\ev} + \tau B_1) + t_2 (\omega_{\od} + \tau B_2)
\]
Therefore, the $\Kee[\tau]$-module of NISes on~$\fg_A(\tau)$ is of rank  $1|1$.

B) Note that the bracket in~$\fg_B(\xi)$ is nonlinear with respect to~$\xi$. There are two NISes on $\fg_B(\xi)$ when $\xi\neq1$:

1) An even NIS~$\omega_e$ for which (and zero on all other pairs of Chevalley basis)
\[\begin{array}{llll}
  &\omega_e (h_1,h_2) = 1, &{}& \omega_e(x_1,y_1) = 1 + \xi^2,\\
  & \omega_e(x_2,y_2) = 1+\xi^2,&{}& \omega_e(x_3,y_3) = 1+ \xi,\\
  &\omega_e (\Pi h_1,\Pi h_2) = 1+\xi, &{}& \omega_e(\Pi x_1,\Pi y_1) = 1 + \xi^2,\\
  & \omega_e(\Pi x_2,\Pi y_2) = 1+\xi^2,&{}&\omega_e(\Pi x_3,\Pi y_3) = 1.
\end{array}
\]

2) An odd NIS~$\omega_o$ for which (and zero on all other pairs of Chevalley basis)
\[\begin{array}{llll}
  &\omega_o (h_1,\Pi h_2) = 1, &{}& \omega_o (h_2,\Pi h_1) = 1+\xi,\\
  &\omega_o(x_1,\Pi y_1) = 1 + \xi^2, &{}& \omega_o(y_1,\Pi x_1) = 1 + \xi,\\
  &\omega_o(x_2,\Pi y_2) = 1+\xi^2, &{}& \omega_o(y_2,\Pi x_2) = 1+\xi,\\
  &\omega_o(x_3,\Pi y_3) = 1, &{}& \omega_o(y_3,\Pi x_3) = 1.
\end{array}
\]

It is easy to see that $\omega_e$ is a~deformation of $\omega_{\ev}$, and $\omega_o$ is a~deformation
of $\omega_{\od}$.

There is no NIS on $\fg_B(\xi)$ for $\xi = 1$:

Thus, the space of NISes on~$\fg_B(\xi)$ is $(1|1)$-dimensional if $\xi \neq 0$ and $0$-dimensional if
$\xi = 1$.
\end{proof}

It now follows that the deform~$\fg_B(1)$ is a~\textit{true} deform, i.e., it is not isomorphic
to~$\fg=\fq(\fsl(3))$, since it has no NISes.

\ssbegin[On $\fk(1; n\vert 1)$]{Lemma}[On $\fk(1; n\vert 1)$]
The Lie superalgebra $\fk(1; n|1)$ and its $(n-2)$-parametric family of even deforms described
  in~\cite[Theorem 6.2]{KL} have no NIS.
\end{Lemma}

\begin{proof}
The even part of the Lie superalgebra $\fk(1;\underline{n}|1)$ and of its $(n-2)$-parameter family of even
deforms is a~solvable Lie algebra~\cite[Corollary~6.4]{KL}. For these Lie superalgebras, we have
$[\fg,\fg]\neq\fg$, and hence no NIS due to Lemma~\ref{1}.
\end{proof}

\section{Remarks}\label{sec:remarks}

\ssec{Multiple NISes over non-closed fields}\label{nonClnis}
NISes on the simple finite-dimensional (associative or Lie) algebra $A$ over algebraically non-closed fields are described in \cite[pp. 30--31, Exercise 15(b)]{Kap} in terms of the centroid  of $A$.  (Recall that the \textit{centroid} of the algebra $A$ is the set of all linear transformations of $A$ that commute with all left and right multiplications. The centroid is an algebra containing the identity linear transformation. If $A$ has a~unit element, the centroid coincides with the elements that commute and associate with everything.) In characteristic 0, Waterhouse considered a~particular case of Kaplansky's exercise: described  NISes on the simple finite-dimensional Lie algebra; they all come from the Killing form over the centroid of the algebra, see \cite{W}. We think that the description below elucidates the above-cited results \cite{Kap, W}, but is more graphic. Superization is immediate; proof is an easy exercise.

\sssec{On NISes over algebraically non-closed fields} \label{nonCl} \textup{1)} Let $\fh$ be a~simple finite-dimensional Lie (super)algebra over $\Ree$ with a~NIS $b$, let ${\fg:=\fh\otimes_\Ree \Cee}$. Clearly,
\[
B(x+iy, z+it):= b(x, z)-b(y,t)+i(b(y,z)+b(x,t))\text{~~for any $x,y,z,t\in\fh$}
\]
is a~NIS on $\fg$.
Then, the realification of $\fg$, i.e., $\fg$ considered as real, 
has at least two  linearly independent NISes
the real and the imaginary parts of  $B$.

 \textup{2)}  More generally, let $\fg$ be a~simple finite-dimensional Lie (super)algebra with a~NIS $B$ over a~field $\Kee$ which is not algebraically closed. For simplicity, let $p\neq 2$.

Let $P \in \Kee[x]$, where $x$ is even, be an irreducible polynomial of degree $d>1$; set ${A := \Kee[x]/(P)}$. The algebra $A$ has no non-trivial ideals and   $A \otimes \fg$ is a~simple Lie (super)algebra usually, but not always. For example, let $\fg$ and $A \otimes \fg$ be simple. Then, $A \otimes (A \otimes \fg)$ is not simple because $A \otimes A$ is not simple.

Let $\varphi: A \tto \Kee$ be a~linear map, not identically equal to zero. Then, the bilinear form on $A \otimes \fg$ given by the formula 
 \[
 B_\varphi(a_1\otimes g_1, a_2\otimes g_2) := \varphi(a_1a_2)B(g_1, g_2)\text{~~for any $a_1,a_2\in A$ and $g_1,g_2\in\fg$}
 \]
 is a~NIS. In this way, we get a~$d$-dimensional space of NISes (of the same parity as $B$) on $A \otimes \fg$ considered as a~Lie (super)algebra over $\Kee$.

\ssec{A Lie superalgebra structure on $L\otimes A$ where $L$ is a~restricted Lie superalgebra and $A$ is a
  non-supercommutative associate superalgebra}\label{sqA}
A more general, but non-super, setting is discussed in \cite{Z}. In this Subsection, $A$ is not necessarily
finite-dimensional.

It is well-known that for any Lie algebra $L$ and any commutative algebra
$A$, one can introduce a~Lie algebra structure on $L\otimes A$ by setting
\[
{}[l_1\otimes a_1, l_2\otimes a_2] := [l_1, l_2]\otimes a_1a_2\text{~~for any $l_i\in L$ and $a_i\in A$.}
\]
The commutativity of $A$ is required for this bracket to be antisymmetric.

Analogously, if $p\neq 2$, $L$ is any Lie superalgebra, and $A$ is any supercommutative
superalgebra, then one can introduce a~Lie superalgebra structure on $L\otimes A$ by setting
\[
[l_1\otimes a_1, l_2\otimes a_2] := (-1)^{p(l_2)p(a_1)}[l_1, l_2]\otimes a_1a_2.
\]
Again, the supercommutativity of $A$ is required for this bracket to be super antisymmetric.

If $p=2$, then the bracket on the Lie superalgebra is again just antisymmetric, so one could assume that the above definition would work even if $A$ is just commutative, not supercommutative (note that if $p=2$, then a~supercommutative superalgebra is commutative as well). But if $p=2$, one has to define the squaring on $(L\otimes A)_\od$ separately from the bracket. For an arbitrary Lie superalgebra $L$, this can be done only if $A$ is a~supercommutative 
superalgebra, and the definition is as follows:

$(l \otimes a)^2 = 0$ for $l\in L_\ev$, $a\in A_\od$;

$(l \otimes a)^2 = l^2\otimes a^2$ for $l\in L_\od$, $a\in A_\ev$;

\[
\left(\sum\limits_{1\leq i\leq n} l_i \otimes a_i\right)^2 = \sum\limits_{1\leq i\leq n} (l_i \otimes a_i)^2 + \sum\limits_{1\leq i<j\leq n} (-1)^{p(l_j)p(a_i)}[l_i, l_j]\otimes a_ia_j
\]
for any homogenous $l_i\in L$ and $a_i\in A$ such that $p(l_i)+p(a_i) = \od$ for all $i$. (The sign $(-1)^{p(l_j)p(a_i)}$ is not needed if $p=2$, but it is introduced here so that the definition would work for other characteristics as well.)

Observe that if $A$ is just commutative, not supercommutative, one can not set
\[
(l \otimes a)^2:=l^2\otimes a^2,
\]
because for any even~$l$ the square is not defined. If $A$ is supercommutative, then for $a$ odd,  $a^2=0$, and everything is OK.

However, if $p=2$ and $L$ is a~Lie superalgebra with a~$2|4$-structure as defined in Subsection~\ref{SSp2pStr} (in particular, if $L$ is a~restricted Lie algebra with a~$2$-structure $l\longmapsto l^{[2]}$ for any $l\in L$), then one can introduce a~Lie superalgebra structure on $L\otimes A$ even if $A$ is a~commutative
superalgebra, as follows:

$(l \otimes a)^2 = l^{[2]}\otimes a^2$ for any $l\in L_\ev$ and $a\in A_\od$;

$(l \otimes a)^2 = l^2\otimes a^2$ for any $l\in L_\od$ and $a\in A_\ev$;

\[
\left(\sum\limits_{1\leq i\leq n} l_i \otimes a_i\right)^2 = \sum\limits_{1\leq i\leq n} (l_i \otimes a_i)^2 + \sum\limits_{1\leq i<j\leq n} [l_i, l_j]\otimes a_ia_j
\]
for any homogenous $l_i\in L$ and $a_i\in A$ such that $p(l_i)+p(a_i) = \od$ for all $i$. (No signs here because $2|4$-structure exists only if $p= 2$.)

Observe that NISes on Lie algebras of the form $L \otimes A$, where $L$ is a~finite-dimensional Lie algebra and $A$ is an associative algebra, are described in \cite[Theorem 4.1]{Z1} and in a~more general situation, where $A$ does not necessarily contain 1, in \cite[Theorem 2]{Z2}.

\end{document}